\documentclass[11pt]{amsart}
\usepackage{amsfonts,amssymb,stmaryrd,amscd,amsmath,latexsym,amsbsy}

\theoremstyle{plain}
\newtheorem{thm}{Theorem}[section]

\newtheorem{prop}{Proposition}[section]

\newtheorem{lem}{Lemma}[section]
\newtheorem{cor}{Corollary}[section]

\newtheorem*{sol}{Solution}

\theoremstyle{definition}
\newtheorem*{defi}{Definition}

\theoremstyle{remark}
\newtheorem{rema}{Remark}[section]

\newcommand{\CC}{\mathbb C}
\newcommand{\ZZ}{\mathbb Z}
\newcommand{\solu}[1]{\begin{sol}{\bf (\ref{#1})}}

\title{ Double affine Hecke algebras and Calogero-Moser spaces}
\author{ Alexei Oblomkov}
\address{Department of Mathematics, MIT, 77, Massachusetts Ave., Cambridge,
MA 02139, USA.}
\thanks{This work was partially supported by the NSF grant DMS-9988796.}
\email{oblomkov@math.mit.edu}

\begin{document}
\begin{abstract}In this  paper  we prove that the spherical subalgebra
$eH_{1,\tau}e$ of the double affine Hecke algebra $H_{1,\tau}$ is
an integral Cohen-Macaulay algebra isomorphic to the center $Z$ of
$H_{1,\tau}$, and $H_{1,\tau}e$ is a Cohen-Macaulay
$eH_{1,\tau}e$-module with the property
$H_{1,\tau}=End_{eH_{1,\tau}e}(H_{1,\tau}e)$. In the case of the
root system $A_{n-1}$ the variety $Spec(Z)$ is smooth and
coincides with the completion of the configuration space of the
relativistic analog of the trigomonetric Calogero-Moser system.
This implies the result of Cherednik that the module $eH_{1,\tau}$
is projective and all irreducible finite dimensional
representations of $H_{1,\tau}$ are regular representation of the
finite Hecke algebra.
\end{abstract}

\maketitle

\section*{Introduction}
Ivan Cherednik in his pioneering  paper \cite{Ch3} introduced
the double affine Hecke algebras. These algebras play a
crucial role in the proof of Macdonald Conjectures \cite{Ch4} and
  are a  natural generalization of  affine Hecke algebras, which are
an object of  great importance in  representation theory.

In the paper \cite{EG} Pavel Etingof and Victor Ginzburg studied
the rational degeneration of a double affine Hecke algebra. They
discovered that in the case when this algebra has a nontrivial
center,  the spectrum of the center is isomorphic to the so called
Calogero-Moser space, and this isomorphism respects the Poisson
structure. The   Calogero-Moser space first appeared in \cite{KKS}
as a completed configuration space for the Calogero-Moser
integrable system.
  Recently  attention to this object was aroused by the paper
\cite{W}.

 The isomorphism between the spectrum of
the center of the degenerate double affine Hecke algebra and the
Calogero-Moser space gives an interpretation of the degenerate
double affine Hecke algebra as  an Azumaya algebra in the case
when the Calogero-Moser space is smooth.

In the present paper we study the double affine Hecke  algebra $H$
with $q=1$. In this case the algebra has a nontrivial center. We
establish a Poisson isomorphism  between the spectrum of the
center $Z(H)$ of $H$ and a relativistic analog of the
Calogero-Moser space in the case of the root system $A_{n-1}$. The
relativistic analog of the Calogero-Moser space is a completed
configuration space for the so called Ruijsenaars-Shneider (or
briefly RS) integrable system \cite{RS}.

For the general algebra $H$ (with $q=1$) we prove that the ring
$Z(H)$ has no zero divisors, and that it is a normal,
Cohen-Macaulay ring
 isomorphic to the spherical subalgebra $eHe$ (where $e$ is the
symmetrizer in the finite Hecke algebra). We also prove the
equality $H=End_{eHe}(He)$, which allows us,  in the case of the
root system $A_{n-1}$, to interpret $H$ as an Azumaya algebra.

The techniques of the paper work also in the degenerate case and
furnish a simpler proof of the results of
 \cite{EG}. Furthermore, there
exists an intermediate degeneration of the double affine Hecke
algebra which lies between the double affine Hecke algebra and the
rational degeneration of this algebra. We call this algebra the
trigonometric degeneration of the double affine Hecke algebra. The
corresponding degeneration of the Calogero-Moser space yields the
configuration space for the trigonometric Calogero-Moser system
(sometimes this space is called the trigonometric Calogero-Moser
space). The results of the paper hold for this intermediate
degeneration and  are given in the last section.

{\bf Acknowledgments.} I'd like to thank my adviser Pavel Etingof
for  posing  the problem,  multiple explanations and help with
proving  many statements from the paper. I am especially grateful
for the  proof of  Theorem 2.
 I am very grateful to Ivan Cherednik who
helped me with the proof of  Lemma~\ref{lcom}. I am  also grateful
to Victor Ostrik for very useful consultations about Hecke
algebras and to  Yuri Berest for his consultations on the
Calogero-Moser space.

\section{Definitions}
\subsection{Definition of the  double affine Hecke algebra
corresponding to $GL(n,\CC)$} We denote this algebra by the symbol
$H_{q,t}$. It is generated by the elements $T_i$, $1\le i\le n-1$,
$\pi$, $X_i^{\pm 1}$, $1\le i\le n$ with relations
\begin{gather}
X_iX_j=X_jX_i,\quad (1\le i,j\le n),\label{dbre1}\\ T_i X_i
T_i=X_{i+1},\quad (1\le i<n),\\ T_iX_j=X_jT_i,\mbox{ if } j-i\ne
0,1\\ [T_i,T_j]=0,\mbox{ if } |i-j|>1\\ T_i T_{i+1} T_i=T_{i+1}T_i
T_{i+1}, (1\le i<n),\\ \pi X_i=X_{i+1}\pi \quad (1\le i\le
n-1),\quad \pi X_n=q^{-1}X_1\pi,\\ \pi T_i=T_{i+1}\pi,\quad \pi^n
T_j=T_j\pi^n,\quad (1\le i<n-1, 1\le j<n)
,\label{dbre2}\\(T_i-\tau)(T_i+\tau^{-1})=0,\quad (1\le i\le
n).\label{T^2}
\end{gather}

\begin{rema}
To identify this definition with the standard definition from the
papers of Cherednik one should replace $\tau$ by $t^{\frac12}$ and
$q$ by $q^{\frac12}$. Also, some definitions use the element
$T_0=\pi T_{n-1}\pi^{-1}$.
\end{rema}
\begin{rema} The double affine Hecke algebra corresponding to $SL(n,\CC)$
is a quotient of the subalgebra of $H_{q,\tau}$ generated by
$X_i/X_{i+1}$, $T_i$, $\pi$, $1\le i\le n-1$, by one extra
relation: $$\pi^n=1.$$
\end{rema}

\subsection{Definition of the Calogero-Moser space}\label{CMsec}
Let $E$ be an $n$-dimensional vector space (over $\CC$). We denote
by the symbol $CM'_\tau$ the subset of  $GL(E)\times GL(E) \times
E\times  E^*$ consisting of the elements $(X,Y,U,V)$ satisfying
the equation
 \begin{equation}\label{CMeq}X^{-1}Y^{-1}XY\tau-\tau^{-1}=U\otimes V.
 \end{equation}
Obviously it is an affine variety.

  The group
$GL(n,\mathbb C)=GL(E)$ acts on it by  conjugation:
\begin{equation*}
(X,Y,U,V)\to (gXg^{-1},gYg^{-1},gU,Vg^{-1}),\quad g\in GL(E).
\end{equation*}
 Later we will show that this action is free if $\tau$ is not
 a root of unity. So the naive quotient  by the action (i.e. the
spectrum of the ring of $GL(E)$ invariant functions) yields an
affine variety, and the quotient is nonsingular if $CM'_\tau$ is.

\begin{defi}
The quotient of $CM'_\tau$ by the action $GL(E)$ is called the
Calogero-Moser space. We use the notation $CM_\tau$ for this
space. \end{defi}

 Below we always suppose that $\tau$ is not a root
of unity.

\section{Properties of the Calogero-Moser space}
The goal of this section is to prove that $CM_\tau$ is a smooth
irreducible algebraic variety of  dimension $2n$. We also
introduce coordinates on  its dense subset. The methods of this
section are analogous to the ones from the paper \cite{W}. In
principal smoothness of $CM_\tau$ follows from the results of the
paper \cite{FR}, the authors of \cite{FR} use the moduli space of the
vector bundles on the punctured torus. For  convenience of reader
we give a direct elementary proof.

\subsection{Smoothness of the Calogero-Moser space}

First we prove a simple  lemma on which all  the following
statements are based.

\begin{lem}If $(X,Y,U,V)\in CM'_\tau$ and $[R,X]=[R,Y]=0$,
$R\in\mathfrak{gl}(E)$ then $R=\lambda Id$ for some $\lambda\in
\mathbb C$.
\end{lem}
\begin{proof}
Let $W\subset E$ be a nonzero subspace which is invariant under
the action of $X$, $Y$ and $R$. We denote by $\bar{X}$ and
$\bar{Y}$ the restriction of the operators $X$, $Y$ to this
subspace. It follows from  equation (\ref{CMeq}) that there are
two possibilities.

In the first case $W\subset V^\perp$, where  $V^\perp$ is the
notation for the annihilator. In this case (\ref{CMeq}) implies
 $$
 \bar{X}^{-1}\bar{Y}^{-1}\bar{X}\bar{Y}=\tau^{-2}Id.
 $$
 But the determinant of LHS is equal to $1$, hence we get a
 contradiction.

 In the second case $W\nsubseteq V^\perp,$ $U\in W$.
  In this case (\ref{CMeq}) implies
 $$
 \bar{X}^{-1}\bar{Y}^{-1}\bar{X}\bar{Y}-U\bar{V}=\tau^{-2}Id,
 $$
 where $0\ne \bar{V}$ is the restriction of $V$ to the
 subspace $W$. Since $det(\bar{X}\bar{Y}\bar{X}^{-1}\bar{Y}^{-1})=1$,
the last equation implies that there is a basis in
 $W$ in which $\bar{X}^{-1}\bar{Y}^{-1}\bar{X}\bar{Y}$ is diagonal
 with the spectrum $\tau^{-2},\tau^{-2},\dots,\tau^{-2},\tau^{2k}$
 where $k=\dim W$. But we know from  equation (\ref{CMeq}) that
 the spectrum of $X^{-1}Y^{-1}XY$ is equal to $\tau^{-2},\tau^{-2},
 \dots,\tau^{2n}$. Thus we get $W=E$.

 The fact that the only common nonzero invariant subspace of $X,Y$ and $R$ is
 the whole $E$  immediately implies the statement of
 the lemma. Indeed, let $\lambda$ be an eigenvalue of $R$, then
 the corresponding eigenspace $W_\lambda$ is invariant under the
 action of $X$ and $Y$, hence it coincides with $E$.
\end{proof}

\begin{cor} The action of $GL(E)$ on $CM'_{\tau}$ is
free.
\end{cor}

\begin{lem} $CM'_\tau$ is smooth.
\end{lem}
\begin{proof}
Let us introduce the  map $\Psi$: $GL(E)\times GL(E)\times E\times
E^* \to \mathfrak{gl}(E)$: $$\Psi(X,Y,U,V)=X^{-1}Y^{-1}XY-U\otimes
V.$$

It is enough to show that $d\Psi$ is epimorphic at a point
$(X,Y,U,V)\in CM'_\tau$. Let $x,y\in\mathfrak{gl}(E)$, $u\in
E,v\in E^*$ and $X(t)=Xe^{xt}$, $Y(t)=Ye^{yt}$, $U(t)=U+tu$,
$V(t)=V+tv$. Then
\begin{multline*} d\Psi_{(X,Y,U,V)}(x,y,u,v)=
\frac{d}{dt}(\Psi(X(t),Y(t),U(t),V(t))|_{t=0}=\\
-xX^{-1}Y^{-1}XY+X^{-1}Y^{-1}XxY-X^{-1}yY^{-1}XY+
X^{-1}Y^{-1}XYy-U\otimes v-u\otimes V.\end{multline*}

If $d\Psi$ is not an epimorphism, then there exists $0\ne
R\in\mathfrak{gl}(E)$ such that
$$tr(d\Psi_{(X,Y,U,V)}(x,y,u,v)R)=0$$ for all
$x,y\in\mathfrak{gl}(E)$, $u\in E,v\in E^*$. Using the cyclic
invariance of the trace, we can rewrite the last condition in the
form: \begin{multline*} tr(x(YRX^{-1}Y^{-1}X-X^{-1}Y^{-1}XYR))+\\
tr(y(RX^{-1}Y^{-1}XY-Y^{-1}XYRX^{-1}))-
v(RU)-VR(u)=0.\end{multline*} As the bilinear form $tr(xy)$ is
nondegenerate, the last equation implies
\begin{gather}
YRX^{-1}Y^{-1}X-X^{-1}Y^{-1}XYR=0,\label{def1}\\
RX^{-1}Y^{-1}XY-Y^{-1}XYRX^{-1}=0,\label{def2}\\ RU=0,\quad
VR=0.\label{def3}
\end{gather}
These equations together with equation (\ref{CMeq}) imply
$[R,X]=[R,Y]=0$. Indeed, let us derive the first equation.

Multiplying on the right formula (\ref{CMeq}) by $R$ we get
\begin{equation}\label{multR}
X^{-1}Y^{-1}XYR=\tau^{-2}R.
\end{equation}
Hence
\begin{multline*}
\tau^{-2}XRX^{-1}=Y^{-1}XYRX^{-1}=RX^{-1}Y^{-1}XY=\\R(\tau^{-1}U\otimes
V+\tau^{-2} Id)=\tau^{-2}Id,
\end{multline*}
here the first equation uses (\ref{multR}), second (\ref{def2}),
third (\ref{CMeq}) and fourth (\ref{def3}).

By the previous lemma $R=\lambda Id$ and finally from (\ref{def3})
we get $R=0$.
\end{proof}
\begin{cor}\label{sm2n} $CM_\tau$ is smooth algebraic variety, and all its
irreducible components have  dimension $2n$.
\end{cor}

\subsection{Local coordinates on $CM_\tau$ }\label{coor}
It is easy to see that  matrices $X,Y\in\mathfrak{gl}(n,\mathbb
C)$,
\begin{gather}\label{repf}
X=diag(\lambda_1,\dots,\lambda_n), \\ \label{repm}Y_{ii}=q_i,\quad
i=1,\dots,n,\\ Y_{ij}=\frac{(\tau-\tau^{-1})q_i\lambda_j}
{(\tau\lambda_i-\tau^{-1}\lambda_j)},\quad 1\le i\ne j\le
n,\label{repl}
\end{gather}
satisfy the equation \begin{equation} \label{eqXY-YX} rk(\tau
XY-\tau^{-1}YX)= 1,\end{equation} for all $\lambda\in(\mathbb
C^*)^n\setminus D_\tau$, $q\in(\mathbb{C}^*)^n$ where $$D_\tau=
\{\lambda|\delta_\tau(\lambda)=\prod_{i\ne
j}(\tau\lambda_i-\tau^{-1}\lambda_j)=0\}.$$

There is a well known formula: if $M=(M_{ij})$, where
$M_{ij}=(\lambda_i-\mu_j)^{-1},$ $ 1\le i,j\le n$, then
\begin{equation*}
\det(M)=\frac{\prod_{i<j}(\lambda_i-\lambda_j)(\mu_j-\mu_i)}
{\prod_{i,j}(\lambda_i-\mu_j)}.
\end{equation*}
To prove this formula one can proceed by the induction on $n$
using the Gaussian method of calculation of the determinant for
the step of the induction.

Applying the last formula to the matrix $Y$ we see that $\det(Y)$
is nonzero if and only if $\lambda_i\ne\lambda_j$, $i\ne j$.

 Let us denote by $\pi'_{12}$: $CM'_\tau\to GL(E)\times
GL(E)$ the projection on the first two coordinates. The previous
reasoning shows that  $(X,Y)\in \pi_{12}(CM'_\tau)$, for
$\lambda\in(\mathbb C^*)^n\setminus (D_\tau\cup D)$, $q\in(\mathbb
C^*)^n$ where $$D= \{\lambda|\delta(\lambda)=\prod_{i< j}
(\lambda_i-\lambda_j)=0\}.$$

Now we can state
\begin{prop}\label{propcoor} Let $(X,Y,U,V)\in CM'_\tau$ and $X$ be
diagonalizable with the different eigenvalues $\lambda_i$,
$i=1,\dots,n$ such that $\tau\lambda_i\ne \tau^{-1}\lambda_j$.
Then the $GL(n,\mathbb C)$ orbit of $(X,Y,U,V)$ contains a
representative satisfying equations $V=\lambda^t$ and
(\ref{repf}-\ref{repl}) for some $q\in(\mathbb C^*)^n$. Such a
representative is unique up to (simultaneous) permutation of the
parameters $(\lambda_i,q_i)$.
\end{prop}
\begin{proof} Equation (\ref{eqXY-YX}) is equivalent to the system
\begin{equation}\label{orb}
\frac{(\tau\lambda_i-\tau^{-1}\lambda_j)Y_{ij}}{\tau-\tau^{-1}}=p_is_j,\quad
1\le i,j\le n, \end{equation} if
$X=diag(\lambda_1,\dots,\lambda_n).$ If there exists $i$ such that
$s_i=0$ then $Y_{ij}=0$, $j=1,\dots,n$ and $\det(Y)=0$. Thus
we have $s_i\ne 0$. Analogously  we get $p_i\ne 0$.

Let us fix a solution of (\ref{orb}) lying in the $GL(n,E)$ orbit
of $(X,Y,U,V)$.  Putting $q_i=p_is_i/\lambda_i$ we get the desired
representative with $X$ given by  formula (\ref{repf}), $Y$ by
formulas (\ref{repm}),(\ref{repl}) and
$U=(\tau-\tau^{-1})X^{-1}Y^{-1}q$.
\end{proof}

This proposition together with Corollary~\ref{sm2n} implies that
$(\lambda,q)$ are local coordinates on the open subset $\bold{U}
\subset CM_\tau$. In the next section we show that this subset is
dense.

\subsection{Irreducibility of $CM_\tau$}
In this subsection we prove
\begin{prop}\label{irrCM} The variety $CM_\tau$ is irreducible.
\end{prop}

  Let us consider the projection on the first
 component $\pi'_1$: $CM'_\tau\to GL(E)$. After the
taking the quotient  by the action of $GL(E)$ this map becomes a
map $\pi_1$: $CM_\tau\to JNF$, where $JNF$ is a stack but we can
think about it as  the set of Jordan normal forms of matrices (we
do not need the stack structure).

Inside $JNF$ there is an open part $\tilde{U}$ corresponding to
diagonal matrices with eigenvalues $\{\lambda_1,\dots,\lambda_n\}$
such that $\lambda_i\ne\lambda_j$,
$\tau\lambda_i\ne\tau^{-1}\lambda_j$ for $i\ne j$. The subset
$\pi_1^{-1}(\tilde{U})$ was described in the previous section. It
is obviously connected. If we show that $\dim
\pi^{-1}(JNF\setminus \tilde{U})<2n$ then Corollary~\ref{sm2n}
implies the irreducibility.

Let us denote by $J_k(\lambda)$ the Jordan block of  size $k$ with
the eigenvalue $\lambda$ and by the symbol $J_{\vec{k}}(\lambda)$
the matrix $diag(J_{k^1}(\lambda),\dots,J_{k^t}(\lambda))$,
$\vec{k}\in\mathbb N^t$ and $k^i\ge k^{i+1}$, $i=1,\dots,t-1$. Let
us formulate without a proof an elementary statement from  linear
algebra.

\begin{lem} The dimension of  $$Stab(J_{\vec{k}}(\lambda))=
\{X\in GL(n,\mathbb C)|[X,J_{\vec{k}}(\lambda)]=0\}$$ is equal to
$\sum_{1\le i,j\le t} min\{k^i,k^j\}$.
\end{lem}

Let us denote by $J_{\bold{k}}(\lambda)$ the matrix
$$diag(J_{\vec{k_1}}(\lambda),J_{\vec{k_2}}(\lambda\tau^2),\dots,
J_{\vec{k_r}}(\lambda\tau^{2r})),$$ $\vec{k_i}\in\mathbb N^{t_i}$.
We  use notations $|\vec{k}_i|=\sum_{j=1}^{t_i} k_i^j$,
$|\bold{k}|=\sum_{j=1}^r|\vec{k}_j|$.

 Let $\lambda_1,\dots,\lambda_s\in\mathbb C$ be such that
$\lambda_i/\lambda_j\ne \tau^{2c}$, $c\in\mathbb Z$, $|c|\le n$
and
\begin{equation}\label{J}
J=
diag(J_{\bold{k}_1}(\lambda_1),\dots,J_{\bold{k}_s}(\lambda_s)).\end{equation}

 We denote by $\pi'_{34}$: $CM'_\tau\to E\times E^*$ the slightly modified projection on the last two components:
$\pi'_{34}(X,Y,U,V)=(YXU,V)$.
 The fiber of the map $\pi'_{34}$ over the point $(U,V)$ of the subset
$\hat{J}=\pi'_{34}((\pi'_1)^{-1}(J))$ consists of the points
$(J,Y+F,J^{-1}(Y+F)^{-1}U,V)$ where $F$ is an element of the
kernel of the linear map:
\begin{equation*}
S_J(F)=\tau J F-\tau^{-1} F J,\quad F\in\mathfrak{gl}(E),
\end{equation*}
 $Y+F$ is invertible, and $(J,Y,J^{-1}Y^{-1}U,V)\in CM'_\tau$.
Obviously  $(\pi'_{34})^{-1}(U,V)$ is a Zariski open nonempty
subset inside
$ker(S_J)$ hence they have the same dimension.

 First let us study the map $S_J$ in the simple case when in
the equation (\ref{J}) we have $s=1$ and $\bold{k}_1=\bold{k}=
(\vec{k}_1,\dots,\vec{k}_r)$, $\vec{k_i}\in\mathbb N^{d_i}$, $1\le
i\le r$. In this situation we denote by $F^{st}_{ij}\in
Mat(k^i_s,k^j_t)$, $1\le s,t\le r$ the matrix with the entries
$F^{st}_{ij;pq}=F_{p'q'}$, $p'=\sum_{l=1}^{s-1}|\vec{k_l}|+
\sum_{l=1}^{i-1} k_s^l+p$, $q'=
\sum_{l=1}^{t-1}|\vec{k_l}|+\sum_{l=1}^{j-1} k_s^l+ q$. In these
notations the following lemma holds

\begin{lem}\label{kerSJ} Let $J$ be the matrix given by (\ref{J}) with $s=1$ and
$\bold{k}_1=\bold{k}=(\vec{k}_1,\dots,\vec{k}_r)$. Then  $F\in
\ker S_J$ if and only if
\begin{gather}
F^{st}_{ij}=0, \mbox{ if } t-s\ne 1,\\ F^{s,s+1}_{ij}=
(\sum_{l=0}^{k_s^i-1}
c^{s}_{ij;l}J_{k_s^i}^l(0))D^{k_s^i,k_{s+1}^j}_\tau \mbox{ if }
k_s^i\le k_{s+1}^j\label{ker1},\\ F^{s,s+1}_{ij}=
D^{k_s^i,k_{s+1}^j}_\tau(\sum_{l=0}^{k_{s+1}^j-1}
c^{s}_{ij;l}J_{k_{s+1}^j}^l(0)) \mbox{ if } k_s^i>
k_{s+1}^j\label{ker2},
\end{gather}
where $c^{s}_{ij;l}\in \mathbb C$, $J^l_{k_s^i}(0)$ ( and
$J^l_{k_{s+1}^j}(0)$) is the $l$-th power of the Jordan block matrix,
and $D^{k_s^i,k_{s+1}^j}\in Mat(k_s^i,k_{s+1}^j)$ is given by
formula
\begin{gather*}
D^{k_s^i,k_{s+1}^j}_{\tau;pq}=\delta_{p+k_{s+1}^j,q+k_s^i}\tau^{2-2p}
\mbox{ if } k_s^i\le k_{s+1}^j\\
D^{k_s^i,k_{s+1}^j}_{\tau,pq}=\delta_{p,q}\tau^{2-2p} \mbox{ if }
k_s^i> k_{s+1}^j.
\end{gather*}
\end{lem}
\begin{proof}
The system of linear equations $S_J(F)=0$ is equivalent to the
collection of  linear systems:
\begin{equation*}
\tau J_{k_s^i}(\lambda\tau^{2s-2})
F^{st}_{ij}-\tau^{-1}F^{st}_{ij}
J_{k_t^j}(\lambda\tau^{2t-2})=0,\quad 1\le s,t \le r,
\end{equation*}
because $J$ has a block structure. The equations for the entries
of $F^{st}_{ij}$ are  of the simple form:
\begin{equation}\label{blok}
F^{st}_{ij;pq}\lambda(\tau^{2s-1}-\tau^{2t-3})=
\tau(\delta_{p,k^s_i}-1)F^{st}_{ij;p+1,q}-
\tau^{-1}(\delta_{q,1}-1)F^{st}_{ij;p,q-1}.
\end{equation}

First consider the case $t-s\ne 1$. Then
$\tau^{2s-1}-\tau^{2t-3}\ne 0$ and equations (\ref{blok}) express
the entries of the $i$-th diagonal through the entries of the
$(i-1)$-th diagonal. It easy to see that in this case (\ref{blok})
implies $F^{st}_{ij;k_s^i,1}=0$, that is, the first diagonal is
zero. Moving from the left to the right we get that all the
diagonals of $F^{st}_{ij}$ are zero.

If $s+1=t$ then equation (\ref{blok}) is a linear relation between
the neighboring entries on the diagonal. It is easy to
derive equations (\ref{ker1}), (\ref{ker2}) from this fact.

Indeed, let us consider the case $k_s^i\le k_{s+1}^j$. Then
equation (\ref{blok}) for $p=k_i^s$, $1< q\le k_i^s$  says
$F^{s,s+1}_{ij;k_i^s,q-1}=0$. Moving along the diagonal from the
bottom to the top and using equation (\ref{blok}) we get that the
first $k_i^s-1$ diagonals of the matrix $F^{s,s+1}_{ij}$ are zero.
For the rest of the diagonals equation (\ref{blok}) implies
$F^{s,s+1}_{ij;p+1,q+p}=F^{s,s+1}_{ij;1,q}\tau^{-2p}$. Putting
$c^s_{ij;l}=F^{s,s+1}_{ij;1,l+k_{s+1}^j-k_s^i+1}$ we get equation
(\ref{ker1}).
\end{proof}

Obviously $Z\in Im S_J$ if and only if $tr(ZF)=0$ for all $F\in
ker \bar{S}_J$, $\bar{S}_J(F)=\tau^{-1} JF-\tau F J$. The space
$ker \bar{S}_J$ has a  description similar to the one of $ker S_J$
(to get $ker\bar{S}_J$ from $ker S_J$ it is enough to change the
order of the Jordan blocks in $J$) and one can easily derive
\begin{cor}
$Z\in Im S_J$ if and only if following equations hold $$
\sum^{u-1}_{l=0}Z_{ij;k_s^i-l,u-l}^{s,s+1}\tau^{2l}=0, \quad
u=1,\dots,min\{k^i_s,k^j_{s+1}\},$$ where
 $s=1,\dots,r-1$.
\end{cor}

 The lowest nonzero diagonal of a  rank one matrix contains  only one
nonzero entry. As $\hat{J}\subset Im S_J\cap \{\mbox{ matrices of
rank } 1\}$ the following statement holds

\begin{cor}\label{Im} $(U,V)\in \hat{J}=\pi_{34}'((\pi'_1)^{-1}(J))$ if
and
only if  $Z=U\otimes V$ satisfies the equation $$
Z_{ij;pq}^{s,s+1}=0 \mbox{ if } p-q\ge
min\{0,k_s^i-k_{s+1}^j\},\quad s=1,\dots,r-1.$$
\end{cor}

Lemma~\ref{kerSJ} gives us the formula for the dimension of the
kernel $$\dim \ker S_J=\sum_{s=1}^{r-1}\sum_{i,j}
min\{k^i_s,k^j_{s+1}\}.$$

We know that $GL(n,E)$ acts on $CM'_\tau$ freely. Hence if we want
to estimate the dimension of the fiber of $\pi_{34}$ over
$\hat{J}$ we should estimate $\dim Stab(J)-\dim ker S_J$. This
difference is positive:
\begin{lem}\label{ineq} Let $k_s\in \mathbb{N}^{d_s}$, $s=1,\dots,r$,
$k^i_s\ge k^{i+1}_{s}$
then the following inequality holds
\begin{equation*}
\sum_{s=1}^r\sum_{i,j} min\{
k_s^i,k_s^j\}-\sum_{s=1}^{r-1}\sum_{i,j} min\{
k_s^i,k_{s+1}^j\}>0,
\end{equation*}
if there  exists $s$ such that $k_s\ne 0$.
\end{lem}
\begin{proof}
Because of the inequality $k^i_s\ge k^{i+1}_s$ we can rewrite LHS
of the inequality in the form
\begin{gather*}
\sum_{\nu=1}\left(\sum_{s=1}^r (x^{\nu}_s)^2- \sum_{s=1}^{r-1}
x^{\nu}_sx^{\nu}_{s+1}\right), \\  x^\nu_s=\#\{i\in\mathbb
N|k_s^i\ge\nu\}.
\end{gather*}
But the first expression  is a sum of positive definite quadratic
forms. Thus we get the lemma.
\end{proof}

The following statement is crucial for  estimating  of
$\dim(\pi^{-1}_1(JNF\setminus\tilde{U}))$:

\begin{prop} If $J$ is given by (\ref{J}) with $s=1$ and $\bold{k}_1=
\bold{k}=(\vec{k}_1,\dots,\vec{k}_r)$,  then $\dim
\pi_1^{-1}(J)<2n-1$ when either $r>1$ or $k_1^1>1$.
\end{prop}
\begin{proof}
In the case $r>1$ Corollary~\ref{Im} implies that $\dim
\pi_{34}'((\pi'_1)^{-1}(J))\le 2n-1$. The theorem on the dimension
of the fibers and previous reasoning imply: $$
\dim\pi_1^{-1}(J)\le \dim\pi'_{34}({\pi'}_1^{-1}(J))+\dim \ker
S_J- \dim Stab(J).$$ Together with the inequality from
Lemma~\ref{ineq} it proves the statement.

Another case (i.e. $\bold{k}=\vec{k}_1$) is even easier  because
in this case we have
\begin{equation*}
\dim\pi_1^{-1}(J)\le 2n -\dim Stab(J)<2n-1.
\end{equation*}
\end{proof}

The case when in  formula (\ref{J}) $s>1$ can be easily reduced to
the previous case. For that  let us introduce the
embedding $i_l$: $\mathfrak{gl}(|\bold{k}_l|,\mathbb C)\to
\mathfrak{gl}(n,\mathbb C)$ and the projection $pr_l$:
$\mathfrak{gl}(n,\mathbb C)\to \mathfrak{gl}(|\bold{k}_l|,\mathbb
C)$: $i_l(Y)_{p',q'}=Y_{pq}$, $pr_l(Y)_{pq}=Y_{p',q'}$,
$p'=p+\sum_{m=1}^{l-1}|\vec{k}_m|$,
$q'=q+\sum_{m=1}^{l-1}|\vec{k}_m|$, $0\le p,q\le |\vec{k}_l|$, and
$i_l(Y)_{ij}=0$ for the rest of the entries of $i_l(Y)$.

Using arguments analogous to the ones from Lemma~\ref{kerSJ} one
gets
\begin{lem} Let $J$ be given by formula (\ref{J}). Then
\begin{enumerate}
\item $ker S_J=\oplus_{l=1}^s i_l(ker
S_{J_{\bold{k}_l}}(\lambda_i))$
\item for $l=1,\dots,s$, $pr_l(Im S_J)\subset Im S_{J_{\bold{k}_l}}$.
\end{enumerate}
\end{lem}

This lemma immediately implies
\begin{prop}\label{dimfib} Let $J$ be given by formula (\ref{J}) and exists $l$,
$1\le l\le s$ such that $|\bold{k}_l|>1$ then
$\dim\pi_1^{-1}(J)<2n-s$.
\end{prop}

And we eventually achieved the goal of the subsection:

\begin{proof}[Proof of Proposition~\ref{irrCM}]
Indeed Proposition~\ref{dimfib} implies
$\dim\pi_1^{-1}(JNF\setminus\tilde{U})<2n$. Hence by
Corollary~\ref{sm2n} $\pi_1^{-1}(JNF\setminus\tilde{U})$ lies
inside the Zariski closure of $\pi_1^{-1}(\tilde{U})$. But
$\pi_1^{-1}(\tilde{U})$  is irreducible.
\end{proof}

\subsection{The Poisson structure on the CM space}\label{Pois}
In the paper \cite{FR} the Poisson structure on the space
$CM_\tau$ was constructed. This Poisson structure on $CM_\tau$
yields the RS integrable system which is the relativistic analog
of the trigonometric Calogero-Moser system.

On the open part $\bold{U}$ of $CM_\tau$ described in the
subsection~\ref{coor} the Poisson bracket $\{\cdot,\cdot\}_{FR}$
takes the form (see Appendix of \cite{FR} for the proof):
\begin{gather*}
\{\lambda_i,\lambda_j \}_{FR}=0,\quad
\{\lambda_j,q_i\}_{FR}=\lambda_i q_i\delta_{ij},\\
\{q_i,q_j\}_{FR}=
\frac{(\tau-\tau^{-1})^2q_iq_j(\lambda_i+\lambda_j)\lambda_i\lambda_j}
{(\tau\lambda_i-\tau^{-1}\lambda_j)(\tau\lambda_j-\tau^{-1}\lambda_i)
(\lambda_i-\lambda_j)}.
\end{gather*}

\begin{rema} The formulas in \cite{FR} contain the misprint, the
authors lost the factor $(\tau^2-1)^2$ in the expression for $\{
q_i,q_j\}_{FR}$.
\end{rema}

 Using the Hamiltonian reduction on the combinatorial model
of the space of flat connections on the torus without a point the
authors of \cite{FR} prove that the Poisson structure
$\{\cdot,\cdot\}_{FR}$ has a holomorphic extension from $\bold{U}$
to the whole $CM_\tau$, and this Poisson structure is
nondegenerate (i.e. $CM_\tau$ is a symplectic variety). Another
way to see this Poisson structure is to use Quasi-Poisson
reduction \cite{AS}. In this picture the Poisson structure is the
result of the reduction of the natural Quasi-Poisson structure on
the product $GL(n,\CC)\times GL(n,\CC)$ and it is immediate that
this Poisson structure is symplectic.

\section{Finite dimensional representation of $H_{1,\tau}$}

In this subsection we construct a family of  finite dimensional
representations of $H_{1,\tau}$. Later we will show that this
family forms an open dense set inside the space of all finite
dimensional representations. The main tool of this section is the
faithful representation of $H_{1,\tau}$ which is the
quasiclassical limit of the standard realization of $H_{q,\tau}$
as a subring of the ring of  reflection difference operators
\cite{Ch4}.

\subsection{Limit of the Lusztig-Demazure operators}

Let us introduce the ring $\tilde{R}=\CC[P_1^{\pm 1}, \dots,
P_n^{\pm 1},X_1^{\pm 1},\dots,X_n^{\pm 1}]_{\delta(X)}\# S_n$,
where the subscript $\delta(X)$ means localization by the ideal
generated by $\delta(X)=\prod_{1\le i<j\le n}(X_i-X_j)$ and $\#$
is a notation for the smash product. Let us explain what the smash
product is. For brevity we will use notation $\CC[P^{\pm 1},
X^{\pm 1}]$ for the  ring $\mathbb C[P_1^{\pm 1}, \dots, P_n^{\pm
1},X_1^{\pm 1},\dots,X_n^{\pm 1}]$

An element of the ring $\tilde{R}$ has the form $\sum_{w\in S_n}
F_w(P,X) w$. The group $S_n$ acts on the ring $R=\mathbb C[P^{\pm
1},X^{\pm 1}]_{\delta(X)}$ by the formulas $$
P_i^w=P_{w^{-1}(i)},\quad X_i^w=X_{w^{-1}(i)},$$ and $$
F(P,X)wF'(P,X)w'=F(P,X)(F')^{w}(P,X)ww'.$$

\begin{prop}\cite{Ch4} The following formulas give an injective homomorphism of
$H_{1,\tau}\to \tilde{R}$:
\begin{gather*}
X^\mu\mapsto X^\mu,\\ T_i\mapsto \tau s_i+ \frac{\tau-\tau^{-1}}
{X_{i}/X_{i+1}-1}(s_i-1), \quad i=1,\dots,n-1,\\ \pi\mapsto
P_1^{-1}c,
\end{gather*}
where $s_i=(i,i+1)\in S_n$ is a transposition and $c\in S_n$ is a
cyclic transformation: $c(i)=i+1, i=1,\dots,n-1$, $c(n)=1$.
\end{prop}

The homomorphism from the proposition is a quasiclassical limit of
the Lusztig-Demazure representation \cite{Ch4}. For brevity we
call this homomorphism the Lusztig-Demazure representation.


\begin{rema} Actually the paper \cite{Ch4} contains the proof for
the case $q\ne 1$. The proof in the case $q=1$ can be obtained
from this proof by  mechanical replacement of shifts operators
$\tau(\lambda)$, $\lambda\in \ZZ^n$ by their quasiclassical limits
$P^\lambda$. The reader may do this operation with Lecture 5 from
the exposition \cite{K}.
\end{rema}

\subsection{The representation $V_{\mu,\nu}$}\label{Vml}
Let $(\mu,\nu)\in (\mathbb C^*)^{2n}$ and $\chi_{\mu,\nu}\simeq
\mathbb C$ be  a one dimensional $R$-module (character):
$\chi_{\mu,\nu}(R(P,X))=R(\mu,\nu)$. We can induce a finite
dimensional module $V_{\mu,\nu}$ from this module:
$$V_{\mu,\nu}=\tilde{R}\otimes_{R}\chi_{\mu,\nu}.$$

This module has a $\mathbb C$  basis $w\otimes 1$, $w\in S_n$,
hence $\dim V_{\mu,\nu}=n!$.

\begin{prop} If $\nu_i\ne \nu_j$, $i\ne j$ then
the $H_{1,\tau}$-module $V_{\mu,\nu}$ is irreducible.
\end{prop}
\begin{proof}
The module $V_{\mu,\nu}$ has a natural  $H_{\delta(X)}$ module
structure. From the Lusztig-Demazure representation we see that
$H_{\delta(X)}\simeq \CC[P^{\pm 1}, X^{\pm 1}]_{\delta(X)}\#S_n$.
The group $S_n$ acts freely on the variety $Spec(\CC[P^{\pm 1},
X^{\pm 1}]_{\delta(X)})$ hence the algebra $H_{\delta(X)}$ is
Morita equivalent to the algebra $\CC[P^{\pm 1}, X^{\pm
1}]^{S_n}_{\delta(X)}$. In particular, the module $V_{\mu,\nu}$
corresponds to the one-dimensional representation: $P\mapsto
P(\mu,\nu)$. Thus $V_{\mu,\nu}$ is an irreducible
$H_{\delta(X)}$-module and hence  an irreducible $H$-module.
\end{proof}

\subsection{The action of the finite Hecke algebra}\label{fHecke}

The elements $T_i$, $i=1,\dots,n-1$ generate an algebra of
dimension $n!$ which is called the finite Hecke algebra. We will
denote it by the symbol $A^n_\tau$.

If $\rm{e}$ is the unit in $S_n$ then by the action of elements
$T_i$ we can get from the vector $\rm{e}\otimes 1$ the whole space
$V_{\mu,\nu}$. Hence the map $j$: $A^n_\tau\to V_{\mu,\nu}$,
$j(T_{i_1}\dots T_{i_k})=T_{i_1}\dots T_{i_k} {\rm e}\otimes 1$ is
an isomorphism of (left) $A^n_\tau$ modules.

\begin{defi}We denote the subset of all finite dimensional irreducible
$H_{1,\tau}$-modules which are regular $A^n_\tau$-modules by the
symbol $Irrep^{n!}$.
\end{defi}

 Let us denote the  subset of $Irrep^{n!}$
consisting of $V_{\mu,\nu}$ $\mu,\nu\in (\CC^*)^n$,
$\delta(\nu)\ne 0$ by $\mathcal U$. Later we will show that all
finite dimensional irreducible modules are from $Irrep^{n!}$.

\subsection{The $GL(2,\mathbb Z)$ action on  double affine Hecke algebras}
\label{Four} One of the most important properties of the double
affine Hecke algebra $H_{q,\tau}$ is the existence  of the action
of $GL(2,\mathbb Z)$
 \cite{Ch4}. To explain how this group acts on the double affine Hecke
algebra we need to introduce pairwise commutative  elements
$Y_i\in H_{q,\tau}$:
\begin{equation}\label{Y}
Y_i=T_1\dots T_{n-i}\pi^{-1} T^{-1}_{n-i+1}\dots T^{-1}_{n-1},
\quad i=1,\dots,n-1.
\end{equation}
These elements satisfy the relations
\begin{gather}\label{YT1}
T_i Y_{i+1} T_i= Y_i, \quad (1\le i<n)\\ T_i Y_j=Y_jT_i,\mbox{ if
} j-i\ne 0,1.\label{YT2}
\end{gather}

The group $GL(2,\mathbb Z)$ is generated by the elements:
\begin{equation*}
\varepsilon=\begin{pmatrix}0&-1\\-1&0\end{pmatrix}, \quad
\sigma=\begin{pmatrix}1&1\\0&1\end{pmatrix}.
\end{equation*}

 These generators  act
by the following formulas:
\begin{gather*}
\varepsilon: X_i\mapsto Y_i, Y_i\mapsto X_i, T_i\mapsto
T_i^{-1},\\ \sigma: X_i\mapsto X_i, Y_i\mapsto X_iY_iq^{-1},
T_i\mapsto T_i,
\end{gather*}
where $\varepsilon: H_{q,\tau}\to H_{q^{-1},\tau^{-1}}$, $\sigma:
 H_{q,\tau}\to H_{q,\tau}$. The transformation $\varepsilon$ is called
 the Fourier-Cherednik transform.

 Using these transformations we can construct  new finite
 dimensional representations. Indeed if $\gamma\in GL(2,\mathbb
 Z)$ is such that $\gamma(H_{1,\tau})=H_{1,\tau'}$ and $\phi':
 H_{1,\tau'}\to GL(V'_{\mu,\nu})$ is the corresponding
 representation of $H_{1,\tau'}$ (here $\tau'$ is either $\tau$ or
 $\tau^{-1}$ ) then the map $\phi'\circ\gamma$ is a
 representation of $H_{1,\tau}$. We denote the set of such representations by
 $\gamma(\mathcal U)$.

\section{The map from $Irrep^{n!}$ to $CM_\tau$}\label{constr}
In this section we construct a map $\Phi$: $Irrep^{n!}\to
CM_\tau$. Later we will show that it is an isomorphism.
Constructions of this section generalize constructions of section
11 of \cite{EG}.

\subsection{Construction of the map}

Let us denote by $A^{n-1}_\tau\subset A^n_\tau$ the subalgebra
generated by the elements $T_{2},\dots, T_{n-1}$. It is the finite
Hecke algebra of  rank $n-2$. The element $v$ of an
$A^n_\tau$-module is said to be $A^{n-1}_\tau$-invariant if
$xv=\tau v$ for all $x\in A^{n-1}_{\tau}.$

The $H_{1,\tau}$-module $V\in Irrep^{n!}$  by definition is a
regular $A^n_\tau$-module. Hence the space $V^{A^{n-1}_\tau}$ of
$A^{n-1}_\tau$-invariants has dimension $n$. The relations inside
$H_{1,\tau}$ and (\ref{YT2}) imply that $X_1$ and $Y_1$ commute
with the action of $A^{n-1}_\tau$. Thus if we fix a basis in $V$
we get $X_1|_{V^{A^{n-1}_\tau}}, Y_1|_{V^{A^{n-1}_\tau}}\in
GL(n,\mathbb C)$.  The following statement is a key statement of
the section.

\begin{prop}\label{rk1} Let $V\in Irrep^{n!}$ then the operators
$\bar{X}_1=X_1|_{V^{A^{n-1}_\tau}}$,
$\bar{Y}_1=Y_1|_{V^{A^{n-1}_\tau}}$ satisfy the equation:
\begin{equation*}
rk(\bar{X}_1\bar{Y}_1\bar{X}^{-1}_1\bar{Y}^{-1}_1-\tau^{-2} Id)=1.
\end{equation*}
\end{prop}

Obviously the space $CM_\tau$ is isomorphic to the quotient of the
space of solutions of (\ref{eqXY-YX}) by the action of
$GL(n,\CC)$. Thus the last proposition proves that the map $\Phi$:
$Irrep^{n!}\to CM_\tau$, $\Phi(V)=(\bar{X}_1,\bar{Y}_1)$ is well
defined.

In the rest  of the section we prove Proposition~\ref{rk1}. It is
done in two steps. First we prove
\begin{lem}\label{lcom} The elements $X_1,Y_1\in H_{1,\tau}$ satisfy the
relation
\begin{equation}\label{com}
X_1Y_1X^{-1}_1Y^{-1}_1=T_1T_2\dots T_{n-2}T_{n-1}^2T_{n-2}\dots
T_1.
\end{equation}
\end{lem}

This is done in the next subsection using the geometric
interpretation of the double affine Hecke algebra. The proof of
this lemma was communicated to the author by Ivan Cherednik. The
last step is the analysis of the LHS of (\ref{com}) using the
quasiclassical limit $\tau\to 1$. It is done in the last
subsection.

\subsection{The double affine  braid group}

The double affine Hecke algebra admits a simple topological
interpretation \cite{Ch3}. This construction is especially simple
in the case $q=1$. In this case the algebra $H_{1,\tau}$ is a
quotient of the so called double affine braid group $\mathfrak
{B}_n$.

For  better understanding of this group the reader may have in
mind the picture analogous to the geometric interpretation of the
usual braid group but in the case when the points live on the two
dimensional torus.

In this picture the elements $T_i$, $i=1,\dots,n-1$ correspond to
the paths which permute of the $i$-th and the $i+1$-th nearby
points (just like in the case of the usual braid group) and $X_i$,
$Y_i$ correspond to the paths in which the $i$-th point goes along
the parallel, respectively  the meridian of the torus. In this
geometric setting it is obvious that  formula (\ref{com}) holds in
$\mathfrak{ B}_n$. The double affine Hecke algebra is a quotient
of $\mathfrak{B}_n$ by the relations (\ref{T^2}). Hence
(\ref{com}) holds in $H_{1,\tau}$. Below we give  formal
definitions to justify this reasoning.

Let $$ U=\{z\in\mathbb C^n| z_k-z_l\notin \mathbb Z+ {\rm
i}\mathbb Z, k\ne l \},$$ and $\bar W=S_n\ltimes(\mathbb
Z\oplus\mathbb Z {\rm i})$ acts on $z\in{U}$ by the formula: $$
\bar{w}({z})=w(z+a+{\rm i} b), \quad \bar{w}=w(a+{\rm i}b), \quad
a,b\in \mathbb Z^n,\quad w\in S_n.$$

We fix a point $z^0$ such that its  real and  imaginary part  is
sufficiently small.

\begin{defi} Paths $\gamma\subset U$ joining $z^0$ with points from
$\{ \bar{w}(z^0),\bar{w}\in \bar{W}\}$ modulo  homotopy and the
action of $\bar{W}$ form the double affine braid group $\mathfrak{
B}_n$ with the multiplication  induced by the usual composing
operation for the paths.
\end{defi}

This group is generated by the elements:
\begin{gather*}
T_j=t_j(\psi)=z^0+(exp(\pi {\rm
i}\psi)-1)(z_j-z_{j+1})(e_j-e_{j+1})\\ X_j=x_j(\psi)=z^0+\psi
e_j,\quad Y_j=y_j(\psi)=z^0+\psi e_j {\rm i},\\ \pi=T_1\dots
T_{n-1} Y_1,
\end{gather*}
where $0\le\psi\le 1$.

\begin{prop} \cite{Ch3}
The group $\mathfrak{ B}_n$ is generated by elements $X_i$, $Y_i$,
$i=1,\dots,n$, $T_j$, $j=1,\dots,n-1$ (and $\pi$) with defining
relations (\ref{dbre1}-\ref{dbre2}), (\ref{YT1}-\ref{YT2}) with
$q=1$.
\end{prop}

Thus $H_{1,\tau}$ is a quotient of $\mathfrak{B}_n$ by relations
(\ref{T^2}) and Lemma~\ref{lcom} follows.

\subsection{The spectrum of $Z=T_1\dots T_{n-2}T_{n-1}^2 T_{n-2}\dots T_1$}
For a representation $V$ is from $Irrep^{n!}$ there is an
isomorphism $V\simeq A^n_\tau$ of  left $A^n_\tau$-modules. Hence
the right multiplication on $A^n_\tau$ induces a structure of a
right $A^n_\tau$-module on $V$ and as a consequence on
$V^{A^{n-1}_\tau}$.

The right  $A^{n}_\tau$-module $V^{A^{n-1}_\tau}$ is  a sum of the
$n-1$ dimensional vector representation and one-dimensional
representation because it is true for $\tau=1$. Obviously, the
operator $Z$ (acting by the left multiplication) commutes with the
right action of $A^n_\tau$. Hence by the Schur lemma $Z$ acts by a
constant on $A^n_\tau$-irreducible components of the right
$A^n_\tau$-module $V^{A^{n-1}_\tau}$. That is, there exists a
basis in the module in which $Z$ is diagonal and of the form
$diag(\lambda_1(\tau),\lambda_2(\tau),\dots,\lambda_2(\tau))$.
Thus we only need to calculate $\lambda_1(\tau)$,
$\lambda_2(\tau)$.

The module $V^{A^{n-1}_\tau}$ exists for all $\tau\ne 0$. As the
operator $Z$ is invertible for all nonzero values of $\tau$, we
have $\lambda_1(\tau)\ne 0$, $\lambda_2(\tau)\ne 0$.

The functions $\lambda_i(\tau)$ are single valued. Indeed for
$n=2$ it is obvious. So let us suppose $n>2$, then the eigenvalues
$\lambda_i(\tau)$ have the different multiplicities. Hence the
Galois group of the extension of the field of rational functions
by $\lambda_1(\tau)$, $\lambda_2(\tau)$ is trivial because it
cannot exchange $\lambda_1$ and $\lambda_2$. Thus the functions
$\lambda_i(\tau)$ are rational and we have $\lambda_i(\tau)=C_i
\tau^{k_i}$, $i=1,2$.

When $\tau=1$, the algebra $A^n_\tau$ becomes the group algebra of
$S_n$, and $Z=1$. Thus we have $C_1=C_2=1$.  The calculation of
$k_1,k_2$ uses the quasiclassical limit reasoning.

If $\tau=e^h$ then we can write the expansion of $T_i$ in terms of
$h$ $$ T_i=s_i+h \tilde{s}_i+O(h^2),\quad i=1,\dots,n-1, $$ where
$s_i=(i,i+1)$ is a usual transposition.  Relation (\ref{T^2})
inside $H_{1,\tau}$ implies
\begin{equation*}
s_i\tilde{s}_i+\tilde{s}_i s_i=2s_i, \quad i=1,\dots,n-1.
\end{equation*}

Let us calculate the first nontrivial term $\tilde{Z}$ of the
expansion of $Z=1+h\tilde{Z}+O(h^2)$:
\begin{equation*}
\tilde{Z}=\sum_{i=1}^{n-1} s_1\dots s_{i-1}(\tilde{s}_i
s_i+\tilde{s}_i s_i)s_{i-1}\dots s_1= 2\sum_{i=1}^{n-1} s_{1i},
\end{equation*}
where $s_{1i}=s_1\dots s_{i-1}s_is_{i-1}\dots s_1$ is a
permutation of $1$ and $i$.

The operator $\tilde{Z}/2$ acts on $\mathbb C[S_n]^{S_{n-1}}$ (by
 the left multiplication) and in the basis $e_i=(\sum_{w'\in
S_{n-1}}w')s_{1i}$ it has the matrix $J-Id$, $J_{ij}=1$, $1\le i,j
\le n.$ Hence $Spec(\tilde{Z}/2)=(n-1,-1,\dots,-1)$. On the other
hand $Spec(\tilde{Z})=(k_1,k_2,\dots,k_2)$. Thus $k_1=2(n-1)$,
$k_2=-2$ and we proved Proposition~\ref{rk1}.

\subsection{The map $\Phi$ on the subset $\mathcal{U}\subset Irrep^{n!}$}
\label{PhiU}

It is possible to calculate $\Phi(V_{\mu,\nu})$ explicitly. Indeed
let us fix a basis in $V_{\mu,\nu}^{A^{n-1}_\tau}$:
$e_i=(\sum_{w'\in S_{n-1}}w')s_{1i}$, $i=1,\dots,n$.

\begin{prop} For the matrices of the operators $\bar{X}_1$ and
$\bar{Y}_1$ written in  the basis $e_i$ the following equations
hold
\begin{gather*}
\bar{X}_1=diag(\nu_1,\dots,\nu_n)\\ \bar{Y}_{ii} =\mu_i
\prod_{j\ne i}\frac{(\tau^{-1}\nu_j-\tau\nu_i)}{(\nu_j-\nu_i)},
\quad i=1,\dots,n.
\end{gather*}
\end{prop}
\begin{proof}
The first equation is obvious. The second formula is a result of
direct calculation using formulas (\ref{Y}) for $Y_1$ and explicit
formulas for $T_i$.

Indeed let make this calculation for $i=1$. The expansion of the
product expression for $Y_1$ consists of the terms of the form
$s_{i_1,j_1}\dots s_{i_r,j_r}c^{-1}F(X)P_1$, where $i_l<j_l$,
$j_m<i_{m+1}$, $l=1,\dots,r$, $m=1,\dots,r-1$ and $F\in
\CC[X^{\pm 1}]_{\delta(X)}$. We know that $Y_1e_1$ is a linear
combination of $e_i$, $i=1,\dots,n$. The terms of the expansion of
$Y_1 e_1$ which contribute to the coefficient before $e_1$ satisfy the
equation $s_{i_1,j_1}\dots s_{i_r,j_r}c^{-1}(1)=1$. This is
possible only in the case $r=1$, $i_1=1,j_1=n$. Thus rewriting
$T_i$ in the form: $$ T_i=\frac{(\tau
X_i-\tau^{-1}X_{i+1})}{X_i-X_{i+1}}s_i+\frac{X_{i+1}(\tau^{-1}-\tau)}
{X_i-X_{i+1}},$$ we see that $$ Y_1e_1=
\left(\prod_{i=1}^{n-1}\frac{(\tau
X_i-\tau^{-1}X_{i+1})}{X_i-X_{i+1}}s_i\right) c^{-1}e_1+R,$$ where
$R$ is a linear combination of $e_j$ with $j>1$. This formula
immediately implies the last formula from the proposition for
$i=1$.
\end{proof}

Let $D_\tau$ be a subset of $\mathcal{U}$ consisting of the
representations of the form $V_{\mu,\nu}$ such that
$\delta_\tau(\nu)=\prod_{i,j}(\tau \nu_i-\tau^{-1}\nu_j)=0$.

 It is actually not easy to compute all coefficients $\bar{Y}_1$
using explicit formulas for $Y_1$ and $T_i$ but we do not need
them. Because by proposition~\ref{propcoor},  if the pair $(X,Y)$
satisfies equation (\ref{eqXY-YX}) and $X$ is diagonal with
eigenvalues satisfying the conditions of
Proposition~\ref{propcoor} then the corresponding $GL(E)$-orbit is
uniquely determined by the diagonal elements of $X$ and $Y$
(because the stabilizer of $X$ consists of  diagonal matrices
which do not change diagonal elements of $Y$ and we can extract
$q$ from these elements). This reasoning implies

\begin{cor}\label{isoopen} The map $\Phi$ is an isomorphism on the subset
$\mathcal{U}_0=\mathcal{U}\setminus D_\tau$, and local coordinates
$\lambda,q$ on $CM_\tau$ are expressed through coordinates
$\mu,\nu$ on $\mathcal{U}_0\subset Irrep^{n!}$ by the formulas
\begin{gather*}
\lambda_i=\nu_i,\quad q_i=\mu_i \prod_{j\ne
i}\frac{(\tau^{-1}\nu_j-\tau\nu_i)}{(\nu_j-\nu_i)}.
\end{gather*}
\end{cor}

\section{ Results
on the general double affine Hecke algebra}

Let $R=\{\alpha\}$ be a root system (possibly nonreduced) of type
$A,$ $B,$ $BC,\dots,F,$ $G$, $W$ the Weyl group generated by the
reflections $s_\alpha$, $\alpha\in R$. The extended affine Weyl
group $\tilde{W}$ is a semidirect product $W\ltimes P$, where $P$
is a weight lattice (i.e. $b\in P$ if
$2(b,\alpha)/(\alpha,\alpha)\in \mathbb Z$ for all $\alpha\in R$).

The affine Hecke algebra $\hat{H}_\tau$ is  a deformation of the
group algebra $\CC[\tilde{W}]$ with deformation parameters
$\tau_\alpha$, $\tau_{w(\alpha)}=\tau_\alpha$, $\alpha\in R$,
$w\in W$ (for the exact definition of the affine Hecke algebra see
\cite{EG}). The double affine Hecke algebra $H_{q,\tau}$ is a
nontrivial extension of the affine Hecke algebra $\hat{H}_\tau$ by
the group algebra $\CC[P^{\vee}]$ of the coweight lattice
$P^{\vee}$ ($b\in P^{\vee}$  if $(b,\alpha)\in\ZZ$ for all
$\alpha\in R$). This extension has one parameter $q$ which is the
shift parameter in the Lusztig-Demazure representation of this
algebra. We consider algebras with $q=1$ and we denote them by
$H$. For the exact definition of the double affine Hecke algebra
and formulas for the Lusztig-Demazure representation see the
original paper \cite{Ch3} or survey \cite{K}.

We use the notation $\delta(X)$ for the Weyl denominator for the
root system $R$. By symbol $\CC[X^{\pm 1}]$ we denote the group
algebra of the weight lattice $P$ lying inside the affine Hecke
algebra $\hat{H}_\tau$ and by symbol $\CC[Y^{\pm 1}]$ we denote
group algebra $\CC[P^{\vee}]\subset H$ which extends
$\hat{H}_\tau$.

There is an injective homomorphism $g$: $H\to\CC[P^{\pm 1} ,X^{\pm
1} ]_{\delta(X)}\#W$ via the quasiclassical Lusztig-Demazure
operators $w\mapsto w, X_b\mapsto X^b,
T_{s_\alpha}=T_{\alpha}\mapsto D_\alpha,$ where $\alpha\in R$,
$b\in P.$ The formulas for the embedding are very similar to the
formulas from the previous section. Let $A$ be the corresponding
finite Hecke algebra, and $e$ the symmetrizer in $A$: $$
e=\sum_{w\in W}\tau^{l(w)} T_w/(\sum_{w\in W}\tau^{2l(w)}),$$
where $T_w=T_{i_1}\dots T_{i_{l(w)}}$ if $w=s_{i_1}\dots
s_{i_{l(w)}}$ is a reduced expression for $w$.

In this section we will need the following PBW type result
\begin{prop}\label{PBW}\cite{Ch4} Each element $h\in H$ can be uniquely
presented in the forms: \begin{gather*}
 h=\sum_{w\in W} f_w(X) T_w
g_w(Y),\\ h=\sum_{w\in W} g'_w(Y) T_w f'_w(X).
\end{gather*}
\end{prop}

\subsection{Formulation of the theorem}
The goal of this section is to study the center $Z$ of $H$ and
corresponding scheme $Spec(Z)$. It turns out that $Z$ is
isomorphic to the subalgebra $eHe$ and we can reduce the study of
$Z$ to the study of $eHe$.

 We remind the definition of a
Cohen-Macaulay algebra.
\begin{defi}\cite{CG} A finitely generated commutative $\CC$-algebra
$A$ is called {\it Cohen-Macaulay} if it contains a subalgebra of
the form $\mathcal{ O}(V)$ such that $A$ is a free $\mathcal
O(V)$-module of finite rank, and $V$ is a smooth affine algebraic
variety.
\end{defi}

 For the definition of a Cohen-Macaulay module see \cite{Ser}
(Chapter 4 p. 18). In this section we prove the following

\begin{thm}\label{general} For any double affine Hecke algebra $H$ the
following is
true:
\begin{enumerate}
\item\label{comm} $eHe$ is commutative.
\item\label{Cohen} $M=Spec(eHe)$ is an irreducible  Cohen-Macaulay and normal
variety.

\item\label{CohenM} The right $eHe$ module $He$ is Cohen-Macaulay.

\item\label{eHe-eH}The left action of $H$ on $He$ induces an isomorphism of
algebras $H\simeq End_{eHe}(He)$.

\item\label{Sat}The map $\eta: z\to ze$ is an isomorphism $Z\to eHe$.
Thus, $M=Spec( Z)$.
\end{enumerate}
\end{thm}

We call the isomorphism $\eta$ the Satake isomorphism (by analogy
with \cite{EG}).

\subsection{Proofs of  theorem~\ref{general}}
\begin{lem} $Z$ contains $\CC[X^{\pm 1}]^W\otimes \CC[Y^{\pm 1}]^W$.
\end{lem}
\begin{proof} $\CC[X^{\pm 1}]^W$ clearly lies in the center of
$\CC[P^{\pm 1},X^{\pm 1}]_{\delta^{-1}(X)}\#W$, and therefore in
the center of $H$. The fact that $\CC[Y^{\pm 1}]^W$ is contained
in $Z$ follows from the existence of the Fourier-Cherednik
morphism \cite{Ch4} (i.e. the action of the element $\varepsilon$
of $GL(2,\ZZ)$ which is described in the previous section in the
case of the root system $A_{n-1}$).

Indeed,  the morphism $\varepsilon$ is an isomorphism between the
double affine Hecke algebra $H'$ with parameter $\tau^{-1}$ and
the double affine Hecke algebra $H$. This morphism maps the
subring $\CC[X^{\pm 1}]^W\subset H'$ onto the subring $\CC[Y^{\pm
1}]^W$ of $H$.

  Now
the statement follows from the PBW theorem for $H$.
\end{proof}

\begin{lem}\label{int} $eHe$ is commutative, without zero divisors.
\end{lem}

\begin{proof}
Let us prove that the subalgebra $eH_{\delta(X)}e$ of
$H_{\delta(X)}\simeq\CC[P^{\pm 1},X^{\pm 1}]_{\delta(X)}\# W$ is
commutative and without zero divisors. Obviously it  implies the
statement.

An element $z\in H_{\delta(X)}$ has a unique representation in the
form $z=\sum_{w\in W} Q_w T_w$; that is, $H_{\delta(X)}$ is
isomorphic to $\CC[P^{\pm 1},X^{\pm 1}]_{\delta(X)}\otimes A$ as a
right $A$-module. If $z\in eH_{\delta(X)}e$ then
$zT_{\alpha}=\tau_\alpha z$ for all $\alpha\in R$ because
$eT_{\alpha}=\tau_\alpha e$. Hence $z$ is an $A$-invariant element
of the right $A$-module $eH_{\delta(X)}e\subset\CC[P^{\pm
1},X^{\pm 1}]_{\delta(X)}\otimes A$. As $\CC(P,X)\otimes A$ is a
regular $A$-module (over the field $\CC(P,X)$) $\CC(P,X)\otimes e$
is  a unique copy of the trivial representation. It implies that
$z=Qe$, $Q\in \CC[P^{\pm 1},X^{\pm 1}]_{\delta(X)}$.

Finally for $z=Qe\in eH_{\delta(X)}e$ we have
$(T_\alpha-\tau_\alpha)Qe=0$. The simple calculation using the
explicit expression for $T_\alpha$ yields: $$
(T_\alpha-\tau_\alpha)Qe=P_\alpha(s_\alpha-1)Qe=P_\alpha
(s_\alpha(Q)-Q)e,$$ where $P_\alpha\in \CC[X^{\pm 1}]_{\delta(X)}$
and $\alpha$ is a simple root. This implies $Q\in \CC[P^{\pm
1},X^{\pm 1}]^W_{\delta(X)}e$ and $eH_{\delta(X)}e\simeq\CC[P^{\pm
1},X^{\pm 1}]^W_{\delta(X)}$.
\end{proof}

The algebra $H$ has a natural  $\CC[X^{\pm 1}]^W\otimes \CC[Y^{\pm
1}]^W$-module structure: the element $p\otimes q$ acts on $x\in H$
by the formula $(p\otimes q)x=pxq$.

\begin{lem}\label{proj} $H$ is a projective finitely generated
$\CC[X^{\pm 1}]^W\otimes \CC[Y^{\pm 1}]^W$-module.
\end{lem}
\begin{proof} Let us first show that $\CC[X^{\pm 1}]$ is a projective finitely
generated $\CC[X^{\pm 1}]^W$ module. Finite generation is clear,
since $W$ is a finite group. Also, it is well known that
$\CC[X]^W$ is a polynomial ring (it is generated by the characters
of the fundamental representations of the corresponding simply
connected group). Since $\CC[X^{\pm 1}]$ is a regular ring, by
Serre's theorem (\cite{Ser}, chapter 4, p. 37, proposition 22)
$\CC[X^{\pm 1}]$ must be locally free over $\CC[X^{\pm 1}]^W$ (in
fact, by Steinberg-Pittie \cite{St} theorem it is free, but we
will not use it). For the same reasons $\CC[Y^{\pm 1}]$ is locally
free over $\CC[Y^{\pm 1}]^W$.

Now the claim follows from the PBW factorization from
Proposition~\ref{PBW} $H=\CC[X^{\pm 1}]\otimes A\otimes \CC[Y^{\pm
1}]$.
\end{proof}

\begin{lem}\label{projeH} $He$ and $eHe$ are projective finitely generated
modules over $\CC[X^{\pm 1}]^W\otimes \CC[Y^{\pm 1}]^W$.
\end{lem}
\begin{proof} The finite generation follows from the Hilbert-Noether lemma and
Lemma~\ref{proj}. The projectivity is true because $He$ and $eHe$
are direct summands in $H$.
\end{proof}

\begin{lem}\label{delta}
\begin{enumerate}
\item $H_{\delta(X)}\simeq \CC[X^{\pm 1},P^{\pm 1}]_{\delta(X)}\# W$
\item  The  map $\eta$: $Z(H_{\delta(X)})\to \CC[P^{\pm 1},
X^{\pm 1}]^W_{\delta(X)}e$, induced by multiplication by $e$ is an
isomorphism.
\item The left $H_{\delta(X)}$-action on $H_{\delta(X)}$ induces
the isomorphism $H_{\delta(X)}\simeq
End_{eH_{\delta(X)}e}(H_{\delta(X)})$.
\end{enumerate}
\end{lem}
\begin{proof}
The first and second items of the lemma follow from the
representation of $H$ by the quasiclassical Lusztig-Demazure
operators. The third item is equivalent to the isomorphism
$$\CC[P^{\pm 1},X^{\pm}]_{\delta(X)}\# W\simeq End_{\CC[P^{\pm
1},X^{\pm 1}]^W_{\delta(X)}}(\CC[P^{\pm 1},X^{\pm
1}]_{\delta(X)}).$$ We will proceed analogously to the proof of
theorem 1.5 from \cite{EG}.

If $a$: $\CC[P^{\pm 1},X^{\pm 1}]_{\delta(X)}\to\CC[P^{\pm
1},X^{\pm 1}]_{\delta(X)}$ is $\CC[P^{\pm 1},X^{\pm
1}]^W_{\delta(X)}$-linear then it defines a $\CC(P,X)^W$-linear
map $\CC(P,X)\to  \CC(P,X)$. The isomorphism $\CC(P,X)\# W\simeq
End_{\CC(P,X)^W}(\CC(P,X))$ implies $a=\sum_{w\in W} a_w w$,
$a_w\in \CC(P,X)$. It is clear that the functions $a_w$ are
regular on $(\CC^*)^n\times (\CC^*\setminus D)^n\setminus \Delta$
where $\Delta$ is the subset of the points of
$(\CC^*)^n\times(\CC^*)^n$ with a nontrivial stabilizer in $W$ and
$D=\{X\in\CC^*|\delta(X)=0\}$. But $\Delta\subset D$, hence
$a_w\in \CC[P^{\pm 1},X^{\pm 1}]_{\delta(X)}$.
\end{proof}

\begin{proof}[Proof of Theorem~\ref{general}]
The first item follows from  Lemma~\ref{int}.

 {\bf Proof of (\ref{Cohen}):} $M=Spec(eHe)$ is an irreducible
affine variety by Lemma~\ref{int}.
 The subalgebra $(\CC[X^{\pm 1}]^W\otimes
\CC[Y^{\pm 1}]^W)e$ is polynomial. Hence to prove that $M$ is
Cohen-Macaulay it is sufficient to show that $eHe$ is a locally
free module of  finite rank over its subalgebra $(\CC[X^{\pm
1}]^W\otimes \CC[Y^{\pm 1}]^W)e$. But the module is projective and
finitely generated by Lemma~\ref{projeH}.

 It is easy to see by localizing with respect to $e\delta(X)$ or
$e\delta(Y)$ that $M$ is smooth away from a codimension 2 subset.
Indeed, by the first item of  Lemma~\ref{delta} after localizing
with respect to $e\delta(X)$ the image of $eHe$ under the
injection $g$ becomes $e\CC[X,Y]_{\delta(X)}e\simeq
\CC[X,Y]_{\delta(X)}^We$, which is the ring of regular functions
on a smooth affine variety. The statement for the localization
with respect to $e\delta(Y)$ follows from the existence of the
Fourier-Cherednik transform. But an irreducible
 Cohen-Macaulay variety that is smooth outside  of a
codimension 2 subset
is normal (\cite{CG},2.2).

{\bf Proof of (\ref{CohenM}):} $eHe$ is finitely generated over
$\CC[X^{\pm 1}]^W\otimes\CC[Y^{\pm 1}]^W$. Hence by Theorem 2.1 of
\cite{BBG} $eH$ is Cohen-Macaulay over $eHe$ if and only if it is
Cohen-Macaulay over $\CC[X^{\pm 1}]^W\otimes\CC[Y^{\pm 1}]^W$.

We know that $He\simeq \CC[X^{\pm 1},Y^{\pm 1}]$ as  a $\CC[X^{\pm
1}]^W\otimes\CC[Y^{\pm 1}]^W$-module and $He$ is projective over
$\CC[X^{\pm 1}]^W\otimes\CC[Y^{\pm 1}]^W$. As $\CC[X^{\pm
1}]^W\otimes\CC[Y^{\pm 1}]^W$ is  a polynomial ring,  the module
$\CC[X^{\pm 1},Y^{\pm 1}]$ is Cohen-Macaulay if and only if it is
projective. So Lemma~\ref{projeH} implies the statement.

 {\bf Proof of (\ref{eHe-eH}):}
 We have an obvious homomorphism $f: H\to End_{eHe}He$. It is
clearly injective because it is injective after localization by
the ideal $(\delta(X))$.

Let us denote $End_{eHe}(He)$ by $H'$.
 Regard $H'\supset H$ as
$\CC[X^{\pm 1}]^W\otimes \CC [Y^{\pm 1}]^W$-modules.
  $H'$ is torsion free because $He$ is  a torsion free
$\CC[X^{\pm 1}]^W\otimes\CC[Y^{\pm 1}]^W$-module (by the PBW
theorem). As $He$ is finitely generated over $eHe$, $H'$ is a
finitely generated $\CC[X^{\pm 1}]^W\otimes\CC[Y^{\pm
1}]^W$-module. Also, $H$ is finitely generated projective, and
$H'/H$ is supported in codimension 2. Indeed, the last part of
Lemma~\ref{delta} implies that $H_{\delta(X)}$ is isomorphic to
$H'_{\delta(X)}$ as a $e H_{\delta(X)} e$ module. Similarly, the
module $H_{\delta(Y)}$ is isomorphic to $H'_{\delta(Y)}$ as
 a $eH_{\delta(Y)}e$-module because we can use (the same way as in the proof of
Lemma~\ref{int}) the Fourier-Cherednik transform $\varepsilon$
from subsection~\ref{Four}.

 The module $H'$ represents some class in
$Ext^1(H'/H,H)$, which must be zero since $H'/H$ is finitely
generated and lives in  codimension 2 and $H$ is projective. Thus,
$H'=H\oplus H'/H$ and the summand $H'/H$ is torsion. But $H'$ is a
torsion free $eHe$ module, hence $H'/H=0$ and $H'=H$.

{\bf Proof of (\ref{Sat}):}
 It is clear that $\eta$ is injective, by looking at the Lusztig-Demazure
representation. Indeed the equation $ze=0$ implies $zp=0$ for any
$p\in \CC[X]^W$,  hence by the PBW theorem $z=0$.

It remains to show that $\eta$ is surjective. Since $eHe$ is
commutative, every element $a\in eHe$ defines an endomorphism of
$He$ over $eHe$ (by right multiplication). So by  statement
(\ref{eHe-eH}) $a$ defines an element $z_a\in H$. This element
commutes with $H$. Indeed, the right multiplication by $a$ is an
endomorphism of the right $eHe$-module which commutes with left
multiplication by elements of $H$ hence by the fourth part of the
theorem $[z_a,h]=0$ for all $h\in H$. For any $x\in H$,
$z_axe=xa$, so $xz_ae=xa$, i.e. $x(z_ae-a)=0$. Since $eHe$ has no
zero divisors, we find $\eta(z_a)=a$, as desired.
\end{proof}

\section{The  results in the case of the root system $A_{n-1}$}
In this section $H=H_{1,\tau}$ is the double Hecke algebra
corresponding to $GL(n,\CC)$.

A point $(\mu,\nu)\in (\CC^*)^n\times ((\CC^*)^n\setminus D)$
defines a $\CC[P^{\pm 1},X^{\pm 1}]_{\delta(X)}^{S_n}$-character
$\chi_{(\mu,\nu)}$: $\chi_{(\mu,\nu)}(Q(P,X))=Q(\mu,\nu)$. The
embedding $Z\hookrightarrow Z_{\delta(X)}\simeq\CC[P^{\pm
1},X^{\pm 1}]_{\delta(X)}^{S_n}$
  allows us to
restrict this character to $Z$. We use the same notation for this
character.

\begin{lem} For any point $(\mu,\nu)\in (\CC^*)^n\times ((\CC^*)^n\setminus D)$
we have $$ He\otimes_{eHe}\chi_{(\mu,\nu)}\simeq V_{\mu,\nu}.$$
\end{lem}
\begin{proof}
The $H$-module $V_{\mu,\nu}$ has a natural structure of an
$H_{\delta(X)}$-module. Let us study finite dimensional
irreducible $H_{\delta(X)}$-modules.

 By  Lemma~\ref{delta} the ring $eH_{\delta(X)}e$ is a
regular ring. As the action of $S_n$ on $(\CC^*)^n\times
(\CC^*\setminus D) $ is free, the  ring $\CC[P^{\pm 1},X^{\pm
1}]_{\delta(X)}\simeq H_{\delta(X)} e$ is a projective
$eH_{\delta(X)}e$-module and defines the vector bundle $F$ over $
\CC^n\times (\CC^n\setminus D)=Spec(eH_{\delta(X)}e)$. Hence by
the last item of Lemma~\ref{delta} $H_{\delta(X)}=End(F)$ is an
Azumaya algebra and by the basic property of  Azumaya algebras any
irreducible $H_{\delta(X)}$-module is of the form
$H_{\delta(X)}e\otimes_{eH_{\delta(X)}e} \chi_{(\mu',\nu')}$ for
some point $(\mu',\nu')\in \CC^n\times (\CC^n\setminus D)$.

Obviously any irreducible $H_{\delta(X)}$-module is irreducible as
an $H$-module. Also we have an obvious isomorphism of $H$ modules
$H_{\delta(X)}e\otimes_{eH_{\delta(X)}e}\chi_{(\mu',\nu')}\simeq
He\otimes_{eHe}\chi_{(\mu',\nu')}$. Thus the previous paragraph
implies $V_{\mu,\nu}\simeq He\otimes_{eHe}\chi_{(\mu',\nu')}$.
Comparing the action of the center on the both sides yields the
statement.
\end{proof}

The previous lemma implies that there is a map $\Upsilon$ from the
open part $Spec(Z_{\delta(X)})$ of $Spec(Z)$ to the  space
$CM_\tau$: $\Upsilon(\mu,\nu)=\Phi(V_{\mu,\nu})$, where $\Phi$ is
the map constructed at the section~\ref{constr}. As
$Spec(Z_{\delta(X)})$ is an open dense subset in $Spec(Z)$, we can
define a rational map $\Upsilon$: $Spec(Z)\dasharrow CM_\tau$.

\begin{thm}The map $\Upsilon$: $Spec(Z)\dasharrow CM_\tau$ is a regular
isomorphism of the algebraic varieties. In particular $Spec(Z)$ is
smooth.
\end{thm}
\begin{proof}
The previous lemma and Corollary~\ref{isoopen} imply that
$\Upsilon$ is a regular isomorphism on
$Spec(Z_{\delta(X)\delta_\tau(X)})$. The Fourier-Cherednik
transform from the section~\ref{Four} allows us to state the same
for the open subset $Spec(Z_{\delta(Y)\delta_\tau(Y)})$.

Indeed, the Fourier-Cherednik transform $\varepsilon$ maps the
double affine Hecke algebra $H_{1,\tau}$ to $H_{1,\tau^{-1}}$ and
it induces the map $\varepsilon_{CM}$: $CM_{\tau}\to
CM_{\tau^{-1}}$,
$\varepsilon_{CM}(X,Y,U,V)=(Y,X,-Y^{-1}X^{-1}YXU,V)$. By the
construction we have
$\varepsilon_{CM}\circ\Upsilon=\Upsilon\circ\varepsilon$. Thus the
restriction of the morphism
$\varepsilon_{CM}^{-1}\circ\Upsilon\circ\varepsilon$ to
$Spec(Z_{\delta(Y)\delta_\tau(Y)})$ is a regular isomorphism.

 Now, we know from the Theorem~\ref{general} that $Spec(Z)$ is
normal. As the complement of
$Spec(Z_{\delta(X)\delta_\tau(X)})\cup
Spec(Z_{\delta(Y)\delta_\tau(Y)})$ has codimension $2$ (because
$Spec(Z)$ is irreducible by Theorem~\ref{general}), we can extend
$\Upsilon$ to a regular map on the whole $Spec(Z)$. The extended
map is dominant because by Proposition~\ref{irrCM} the variety
$CM_\tau$ is irreducible.

Thus $\Upsilon$ is a regular birational map which is an
isomorphism outside of the subset of codimension $2$.  But we know
that $CM_\tau$ is smooth and $Spec(Z)$ is normal, hence (by
theorem~5 section~5 of chapter~2 of \cite{Sh}) the map
$\Upsilon^{-1}$ is regular and as a consequence is an isomorphism.

\end{proof}

\begin{cor} $He$ is a projective $eHe$-module.
\end{cor}
\begin{proof} We proved for any $R$ that $He$ is a Cohen-Macaulay module
over $eHe$. Since  $M=Spec(eHe)$ is smooth, the result follows
from corollary 2 from chapter 4 of \cite{Ser}.
\end{proof}

Thus $He$ defines the vector bundle $E$ on $Spec(eHe)$, with
fibers of  dimension $n!$.

\begin{cor}\label{Azum} For the double affine Hecke algebra
$H=H_{1,\tau}$ the following is true: \begin{enumerate}
 \item $H=End
E$ where $E$ is a vector bundle over $Spec(Z)$ i.e. $H$ is an
Azumaya algebra.
\item  Every irreducible
representation of $H$ is of the form $V_z=He\otimes_{eHe}\chi_z$,
 $z\in M=Spec(Z)$.
 \item $V_z$ has dimension $n!$ and is a regular representation of
$A^n_{\tau}$.
\end{enumerate}
\end{cor}
\begin{proof} The first item follows from
Theorem~\ref{general}. The second item is a general property of
Azumaya algebras. The third item follows from the fact that it is
true for the generic point $z\in Spec(Z)$.
\end{proof}

\begin{rema} This corollary was proved in 2000 by Cherednik using
different methods \cite{Chpriv}.
\end{rema}

The ring $Z\simeq eH_{1,\tau}e$ has a natural noncommutative
deformation $e H_{q,\tau}e$. Hence this ring has a natural Poisson
structure $\{\cdot,\cdot\}$. The variety $CM_\tau$ also has a
Poisson structure described in  subsection~\ref{Pois}. It turns
out that the isomorphism $\Phi$ respects these Poisson structures.

\begin{thm}The  isomorphism $\Phi$ is an isomorphism of  Poisson
varieties, that is  the following formula holds $$
\{\cdot,\cdot\}_{FR}=\{\cdot,\cdot\}. $$
\end{thm}
\begin{proof}
It is enough to prove that it is an isomorphism of Poisson
varieties on the open set $\mathcal{U}$. For $q=e^h\ne 1$ we have
an embedding $g_q$; $H_{q,\tau}\to \mathbb{D}_q\# S_n$ via
Lusztig-Demazure reflection difference operators. Here
$\mathbb{D}_q$ is a localization of the Weyl algebra with
generators $X_i^{\pm 1}$, $\hat{P}^{\pm 1}_i$, $i=1,\dots,n$ and
relations: $$ [X_i,X_j]=0,\quad [\hat{P}_i,\hat{P}_j]=0,\quad
X_j\hat{P}_i-q^{\delta_{ij}}\hat{P}_iX_j=0,$$ by the ideal
$(\delta(X))$. When $q=1$, the noncommutative ring $\mathbb{D}_q$
becomes the commutative ring $\CC[P^{\pm 1},X^{\pm
1}]_{\delta(X)}$
 and the corresponding Poisson structure on this ring
 is given by the formulas:
\begin{equation*}
\{ X_i,X_j\}=0,\quad\{ P_i,P_j\}=0,\quad
\{X_i,P_j\}=\delta_{ij}X_iP_j.
\end{equation*}

The $H_{1,\tau}$-module $V_{\mu,\nu}$ has a natural $\CC[P^{\pm 1}
,X^{\pm 1}]_{\delta(X)}\#S_n$ structure. It is easy to see that in
the basis $1\otimes w$, $w\in W$
 operators $P_i$, $X_j$ are diagonal.
 In particular $P_i(1\otimes {\rm e})=\mu_i(1\otimes{\rm e})$ and
 $X_i(1\otimes{\rm e})=\nu_i(1\otimes{\rm e})$, hence we have
 the following Poisson bracket on $\mathcal{U}$:
\begin{equation}\label{mnPoi}
\{\nu_i,\nu_j\}=0,\quad \{\mu_i,\mu_j\}=0\quad
\{\nu_i,\mu_j\}=\delta_{ij}\nu_i\mu_j.
\end{equation}
The  comparison of the formulas for the Poisson bracket on
$\bold{U}\subset CM_\tau$ from  subsection~\ref{Pois} and explicit
formulas for the map $\Phi|_{\mathcal{U}}$ from
subsection~\ref{PhiU} give the formula. Indeed, we can express the
functions $\lambda_i$, $q_k$ through the functions $\mu_s,\nu_t$
and using (\ref{mnPoi}) calculate the Poisson brackets
$\{\lambda_i,\lambda_k\}$, $\{\lambda_i,q_k\}$, $\{q_i,q_k\}$. We
prove a formula for the last bracket:
\begin{multline*}
\{q_i,q_k\}=q_iq_k\left(\nu_k\frac{\partial \ln(q_i)}{\partial
\nu_k}-\nu_i\frac{\partial\ln(q_k)}{\partial\nu_i}\right)=\\
q_iq_k\left(\nu_k\left(-\frac{\tau}{\tau^{-1}\nu_i-\tau\nu_k}+\frac{1}{\nu_i-\nu_k}\right)-
\nu_i\left(-\frac{\tau}{\tau^{-1}\nu_k-\tau\nu_i}+
\frac{1}{\nu_k-\nu_i}\right)\right)=\\
\frac{(\tau^{-1}-\tau)^2q_iq_k(\nu_k+\nu_i)\nu_i\nu_k}{(\nu_i-\nu_k)
(\tau^{-1}\nu_k-\tau\nu_i)(\tau^{-1}\nu_i-\tau\nu_k)}.
\end{multline*}
\end{proof}

\section{The rational and trigonometric cases}
In this section we explain how to degenerate results from the main
body of the paper to obtain an easier proof of the results of
\cite{EG} on the rational double affine Hecke algebra. We also
give the version of the results of the paper for the trigonometric
Hecke algebra and explain how to modify the proof from the paper
for this case.

We give the modifications of the results from the main body of the
text only for the root system $A_{n-1}$ but  similar things can be
done for any root system $R$. Moreover, in the rational case we
can replace the Weyl group $W$ by a finite Coxeter group (see
\cite{EG}). Proofs of these results almost identically repeat
proofs for (nondegenerate) double affine Hecke algebras.

\subsection{Definition of the rational and trigonometric
 double affine Hecke algebras}
Below we give a definition of the rational and trigonometric
double affine Hecke algebra.

\begin{defi}\cite{EG,ChM} The rational double affine Hecke algebra $H^{rat}_{t,c}$ is
 generated
by elements $s_{ij}$, $1\le i\ne j\le n$, $x_i,y_j$, $1\le i,j\le
n$. The elements $s_{ij}$, $1\le i,j\le n$ generate the subalgebra
inside $H^{rat}_{t,c}$  isomorphic to the group algebra of the
symmetric group $S_n$, and $s_{ij}$ corresponds to the
transposition $(ij)$. In addition generators  of $H^{rat}_{t,c}$
satisfy the relations
\begin{gather*}
x_is_{ij}=s_{ij}x_j,\quad y_is_{ij}=s_{ij}y_j,\quad 1\le i,j\le
n,\\ [x_k,s_{ij}]=0,\quad [y_k,s_{ij}]=0,\quad k\notin \{i,j\},
\quad 1\le i,j,k\le n,\\ [y_i,x_j]= cs_{ij},\quad 1\le i\ne j\le
n\\ [x_i,x_j]=0=[y_i,y_j],\quad 1\le i, j\le n,\\
[y_k,x_k]=t-c\sum_{i\ne k}s_{ik},\quad 1\le k\le n.
\end{gather*}
\end{defi}

\begin{defi} The trigonometric double affine Hecke algebra $H^{trig}_{t,c}$ is
generated by elements  $s_{ij}$, $1\le i\ne j\le n$, $X^{\pm 1}
_i,y_j$, $1\le i,j\le n$. The elements $s_{ij}$, $1\le i,j\le n$
generate the subalgebra inside $H^{rat}_{t,c}$  isomorphic to the
group algebra of the symmetric group $S_n$, and $s_{ij}$
corresponds to the transposition $(ij)$.  In addition the
generators of  $H^{trig}_{t,c}$ satisfy the relations
\begin{gather*}
X_is_{ij}=s_{ij}X_j,\quad 1\le i,j\le n,
\\ s_{ij}y_i-y_js_{ij}=c \mbox{ if } j>i,\quad
s_{ij}y_i-y_js_{ij}=-c\mbox{ if } j<i,\\ [X_k,s_{ij}]=0,\quad
[y_k,s_{ij}]=0 \mbox{ if } k\notin \{i,j\}, \quad 1\le i,j,k\le
n,\\ [X_i,X_j]=0= [y_i,y_j],\quad 1\le i, j\le n,\\
X_j^{-1}y_iX_j-y_i=cs_{ij}\mbox{ if } j>i,\quad
X_j^{-1}y_iX_j-y_i=X_iX_j^{-1}cs_{ij}\mbox{ if } j<i,\\
 X_k^{-1}y_kX_k-y_k=t-c(\sum_{i<k
}s_{ik}+\sum_{i>k}X_iX_k^{-1}s_{ik}),\quad 1\le k\le n.
\end{gather*}
\end{defi}

\begin{rema}  Let
$\hat{H}$  be the  $\CC[c,t][[h]]$-algebra topologically generated
(in the $h$-adic topology) by $X_i$, $y_i$, $s_{i,i+1}$ with
$T_i=s_{i,i+1}e^{chs_{i,i+1}}$, $i=1,\dots,n-1$, $Y_i=e^{hy_i}$,
$X_i$, $i=1,\dots,n$ satisfying the relations for the double
affine Hecke algebra $H_{q,\tau}$, $q=e^{th},\tau=e^{ch}$. It
coincides with an appropriate completion of the double affine
Hecke algebra $H_{q,\tau}$,
 in the $h$-adic topology. Moreover, one can show that
$\hat{H}$ is flat over $\CC[[h]]$ and
$\hat{H}/h\hat{H}=H^{trig}_{t,c}$. Analogously,
  if $\hat{H}^{trig}$  is the $\CC[c,t][[h]]$-algebra topologically
generated
  by  by $s_{ij}, y_i, x_j$, $1\le i\le n$ with $s_{ij}, y_i,X_j=e^{hx_j},$
  $i,j=1,\dots,n,$ satisfying the relations for the trigonometric double
  affine Hecke algebra $H^{trig}_{ht,hc}$ then  the algebra
  $\hat{H}^{trig}$ is flat over $\CC[[h]]$ and
  $H^{rat}_{t,c}=\hat{H}^{trig}/h\hat{H}^{trig}$.
\end{rema}

\subsection{Representation  by   Dunkl operators} Let
$\mathcal{D}^{rat}_t$ be the localization
 of the $n$-dimensional Weyl  algebra $\mathcal{A}^{rat}_t$ by the ideal
generated by
$\delta(x)$. The Weyl algebra $\mathcal{A}^{rat}_t$ is generated
by elements $x_i,p_i$, $1\le i\le n$ modulo relations:
\begin{equation*}
[x_i,x_j]=0=[p_i,p_j],\quad [x_i,p_j]=t\delta_{ij},\quad 1\le
i,j\le n.
\end{equation*}

Let us denote by $\mathcal{D}_t^{trig}$ the trigonometric version
of the algebra $\mathcal{D}^{rat}_t$. This algebra is localization
by $(\delta(X))$ of the algebra $\mathcal{A}_t^{trig}$ with
generators $p_i, X^{\pm 1}_i$, $i=1,\dots,n $ modulo relations:
\begin{equation}\label{AA}
[X_i,X_j]=0=[p_i,p_j],\quad [X_i,p_j]=t\delta_{ij}X_i,\quad 1\le
i,j\le n.
\end{equation}
It is easy to see that the ring $\mathcal{A}_t^{trig}$ is
isomorphic to the ring of differential operators on the torus
$(\CC^*)^n$.

\begin{prop}\label{Dunklpr} The homorphisms $g^{rat}$: $H^{rat}_{t,c}
\to  D^{rat}_t \#S_n$, $g^{trig}$: $H^{trig}\to
\mathcal{D}^{trig}_t\#S_n$  defined by the formulas
\begin{gather*}
g^{rat}(y_i)=p_i+c\sum_{j\ne i}\frac1{x_i-x_j}(s_{ij}-1),\\
g^{rat}(x_i)= x_i,  \quad
 g^{rat}(w)=  w,\\
g^{trig}(y_i)= p_i+
c\sum_{j<i}\frac{X_i}{X_i-X_j}(s_{ij}-1)+c\sum_{j>i}
\frac{X_j}{X_i-X_j}(s_{ij}-1),\\ g^{trig}(X_i)=X_i, \quad
  g^{trig}(w)= w,
\end{gather*}
 ($i=1,\dots,n$) is injective.
\end{prop}

This proposition allows  to prove the PBW type result for these
algebras.

\subsection{Calogero-Moser spaces}
In this subsection we give a definition of the Calogero-Moser
space in the rational and trigonometric cases.

Let $CM'_{rat}$ be the subset of
$\mathfrak{gl}(n,\CC)\times\mathfrak{gl}(n,\CC)$ consisting of the
elements $(x,y)$ satisfying the equation $$ rk([x,y]+Id)=1. $$ By
$CM'_{trig}\subset GL(n,\CC)\times \mathfrak{gl(n,\CC)}$ we denote
the subset of pairs $(X,y)$ satisfying: $$ rk(X^{-1}yX-y+Id)=1.$$

The group $GL(n,\CC)$ acts on the spaces $CM'_{rat}$ and
$CM'_{trig}$ by conjugation. This action is free.

\begin{defi}
The quotient of $CM'_{rat}$ ($CM'_{trig}$) by the action of
$GL(n,\CC)$ is called the rational (trigonometric) Calogero-Moser
space. We use the notation $CM_{rat}$ (respectively $CM_{trig}$)
for this space.
\end{defi}

\begin{prop} The rational (trigonometric) Calogero-Moser space
$CM_{rat}$ ($CM_{trig}$) is an irreducible smooth variety of
dimension $2n$.
\end{prop}

For the rational Calogero-Moser space this statement is proved in
section 1 of \cite{W}. The proof in the trigonometric case almost
identically repeats the proof in the rational case.

The  Calogero-Moser spaces $CM_{rat}$ and $CM_{trig}$ are the
configuration spaces for the rational and trigonometric integrable
Calogero-Moser systems. The Poisson structures corresponding to
these systems are the results of the Hamiltonian reduction of the
natural Poisson structures on the spaces
$\mathfrak{gl}(n,\CC)\oplus \mathfrak{gl}^*(n,\CC)$ and
$T^*GL(n,\CC)$ (see \cite{OP}).

\subsection{The main result for the rational and trigonometric
double-affine Hecke algebras}

As we mentioned in the first subsection, the algebras
$H^{rat}_{0,c}$ $H^{trig}_{0,c}$ are in some sense quasiclassical
limits of the double affine Hecke algebra $H_{1,\tau}$. Naturally,
the theorems from the previous section have their rational and
trigonometric analogs:

\begin{thm} Let $H$ be one of three described algebras:
$H_{1,\tau}$, $H^{trig}_{0,c}$, $H^{rat}_{0,c}$, $CM$ is the
corresponding Calogero-Moser space, and $e$ is the symmetrizer (in
the finite Hecke algebra if $H=H_{1,q}$ and in the symmetric group
otherwise). Then the following is true:
\begin{enumerate}
\item The map $h$: $z\to z e$ is an isomorphism between $Z(H)$ and
$eHe$.
\item $Spec(Z(H))$ is an irreducible smooth variety naturally isomorphic to
$CM$.
\item The Poisson structure on $CM$ which comes from the noncommutative
deformation $eH_{q,\tau}e$ ($eH^{trig}_{t,c}e$, $eH^{rat}_{t,c}e$
respectively) of $eHe$ coincides (up to a constant) with the
(Quasi) Poisson structure on $CM$ coming from the (Quasi)
Hamiltonian reduction.

\item The left $eHe$-module $He$ is projective and
$H=End_{eHe}(He)$.
\end{enumerate}
\end{thm}

In particular the algebras $H^{rat}_{0,c}$ and $H^{trig}_{0,c}$
are Azumaya algebras and for these algebras the statement of
Corollary~\ref{Azum} holds with $A^n_\tau$ replaced by $S_n$.

The proof of the theorem in the case $H=H^{rat}_{0,c}$ is
completely parallel  to the case $H=H_{1,\tau}$.

 In the trigonometric case the only difficulty is that the
 group $GL(2,\ZZ)$ does not act on
$H_{0,c}^{trig}$ and we do not have any analog of the
Fourier-Cherednik transform. But instead of the Fourier-Cherednik
transform one can use the faithful representation $\bar{g}^{trig}$
of  $ H_{0,c}^{trig}$. The representation $\bar{g}^{trig}$ is
the "bispectral dual"  to $g^{trig}$; that is, the role of $X_i$,
$1\le i\le n$ is played by $y_i$, $1\le i\le n$.

Let us describe the representation $\bar{g}^{trig}$.  The
homomorphism $\bar{g}^{trig}: H_{t,c}^{trig}\to \CC[P^{\pm
1},y]_{\delta(y)}\# S_n$ is defined by the formulas
\begin{gather*}
s_{i,i+1}\mapsto
\bar{T}_i=s_{i,i+1}+\frac{c}{y_i-y_{i+1}}(s_{i,i+1}-1),\quad 1\le
i\le n-1,
\\ y_i\mapsto y_i,\quad 1\le i\le n,\\
X_i\mapsto \bar{T}_1\dots \bar{T}_{n-i}wP_1\bar{T}_{n-i+1}\dots
\bar{T}_{n-1},\quad 1\le i\le n,
\end{gather*}
where $w\in S_n$, $w(1)=n$, $w(i)=i-1$, $i=2,\dots,n$.

\end{document}